\documentclass[reqno,12pt]{amsart}

\newtheorem{theorem}{Theorem}[section]
\newtheorem{lemma}[theorem]{Lemma}
\newtheorem{corollary}[theorem]{Corollary}

\theoremstyle{definition}

\newtheorem{exercise*}[theorem]{Exercise*}

\theoremstyle{remark}

\newtheorem{remark}[theorem]{Remark}

\makeatletter
\def\dashintindex{\operatorname%
{-\kern-.7em\DOTSI\intop\ilimits@}}%
\def\dashint{\operatorname%
{\,\,\text{\bf--}\kern-.98em\DOTSI\intop\ilimits@\!\!}}
\makeatother

\newcommand\bH{\mathbb{H}}

\newcommand\bR{\mathbb{R}}

\newcommand\cA{\mathcal{A}}

\newcommand\cH{\mathcal{H}}

\newcommand\cL{\mathcal{L}}
\newcommand\cM{\mathcal{M}}

\newcommand{\Div}{{\rm div}\,}
\newcommand{\osc}{{\rm osc}\,}

\newcommand{\cHO}{\overset{\scriptscriptstyle0}%
{\mathcal H}\,\!}
\newcommand{\WO}{\overset{\scriptscriptstyle0}%
{W}\,\!}

 \newcommand{\mysection}[1]{\section{#1}
 \setcounter{equation}{0}}

\begin{document}

\title%
{Parabolic and elliptic equations with VMO coefficients}
\author{N.V. Krylov}
\thanks{The work  was partially supported by
NSF Grant DMS-0140405}
\email{krylov@math.umn.edu}
\address{127 Vincent Hall, University of Minnesota, Minneapolis, MN, 55455}
 
 

\begin{abstract}
An $L_{p}$-theory
 of divergence and non-divergence form
elliptic  and parabolic
equations is presented.
The main coefficients are supposed to belong to the
class $VMO_{x}$,
which, in particular, contains all functions
independent of $x$. Weak uniqueness of the martingale
problem associated with such equations is obtained. 
\end{abstract}

\maketitle

\mysection{Introduction}

The goal of this paper is to expand the $L_{p}$-theory
of parabolic equations to a larger class of operators,
having discontinuous coefficients, than previously  known.
By doing this we also obtain a generalization of a result
of Stroock-Varadhan \cite{SV} about
weak uniqueness of solutions of It\^o equations, which
was our main motivation
(see Remark \ref{remark 7.21.1}). For issues related to
  stochastic processes it is enough to consider
the corresponding PDEs in the whole space, and therefore
we do not consider boundary-value problems.
Uniqueness problem for stochastic equations
is an old one. Recently the interest in solving it
for discontinuous coefficients 
reappeared in connection with
diffusion approximation (see \cite{Kr04}, \cite{KL}).

According to the famous counterexample of 
Nadirashvili there could not exist theory of solvability
of equations with general discontinuous coefficients
even if they are uniformly bounded and equations are
uniformly elliptic. Therefore, much effort was
applied to treat particular cases of discontinuity.
First came equations with piecewise continuous coefficients,
see \cite{Lo} and \cite{Ki}. Then,  truly remarkable
and absolutely unpredictable results
about $W^{2}_{p}$-estimates for elliptic equations with
VMO coefficients appeared in \cite{CFL1}. They were later
developed into existence theory for {\em non-divergence\/}
 form elliptic and parabolic equations in
\cite{CFL2} and \cite{BC}. The results in \cite{BC},
\cite{CFL1}, and \cite{CFL2} are based on deep versions
of the Calder\'on-Zygmund theorem and estimates
of certain commutators.

A different approach 
to {\em divergence\/} form elliptic and parabolic equations
with VMO coefficients is developed in \cite{By1} and \cite{By2}.
These two papers also could  be used as a good source of
further references on the subject of VMO and
equations in divergence and non-divergence forms.
One can also consult papers \cite{HHH} and \cite{Li}
 for various versions and extensions.

In what concerns parabolic equations there is a flaw in
the results in \cite{BC}, \cite{By2},
 and  \cite{HHH}. Namely, these
  do {\em not\/} contain quite classical results
about solvability in Sobolev spaces of equations
whose leading coefficients  depend  only on $t$
and are just {\em measurable\/} functions
(see  \cite{SV} and the references therein).

We correct this flaw and 
treat divergence and non-divergence form
elliptic  and parabolic
 operators,
the main coefficients of which are in VMO.
Actually, as in \cite{By1} and \cite{By2},
 a slightly more general
class of coefficients  is allowed (see Remark \ref{remark
7.23.1}). In contrast with many of the above references  we do
not  consider boundary-value problems for the  reasons
explained in the beginning. This also makes
the presentation clearer and allows us to use
a unified approach to elliptic and parabolic
divergence and non-divergence form equations.
We do not treat the $L_{q}-L_{p}$ theory either
and only mention that parabolic equations 
with mixed norms and coefficients constant in time
are considered in \cite{HHH} and with coefficients,
which are uniformly {\em continuous\/} in $x$
and measurable in $t$,
in \cite{Kr02}.

In a sense our approach is a combination of
the approach in \cite{CFL1}, \cite{CFL2},
\cite{BC} and the one in \cite{By1}, \cite{By2}.
On the one hand, we use pointwise estimates
of the sharp function of second order derivatives, 
on the other hand, we do not use integral representations
of these derivatives to deal with
contributions from ``far away",
but deal with these contributions by splitting the function into
two parts, one of which is ``harmonic", that is
satisfies the homogeneous equation. 

The article is organized as follows. Section
\ref{section 7.22.1} contains our main results.
In Section \ref{section 7.22.2} we discuss
some auxiliary results that are later used
for non-divergence and divergence form equations.
Section \ref{section 7.22.3} is devoted to proving
our main results for non-divergence type equations.
Then comes Section \ref{section 7.22.4} with few
more auxiliary results needed for divergence type
equations and the short final Section \ref{section 7.22.5}
deals with the proofs of our results for such 
equations.

Hongjie Dong and Doyoon Kim kindly showed the author
few errors in the original version of the article
for which the author is sincerely grateful.

\mysection{Main results}
                                      \label{section 7.22.1}

Let $\bR^{d}$ be a $d-$dimensional Euclidean space of points
$x=(x^{1},...,x^{d})$  and
$$
\bR^{d+1}=\{(t,x):t\in\bR,x\in\bR^{d}\}.
$$

We are dealing with four types of operators:
parabolic
$$
Lu(t,x)=u_{t}(t,x)+a^{ij}(t,x)u_{x^{i}x^{j}}(t,x)+
b^{i}(t,x)u_{x^{i}}(t,x)+c(t,x)u(t,x) ,
$$
$$
\cL u(t,x)=u_{t}(t,x)+\big(a^{ij}(t,x)u_{x^{i}}(t,x)+
\hat{b}^{j}(t,x)u(t,x)\big)_{x^{j}} 
$$
$$
+b^{i}(t,x)u_{x^{i}}(t,x)+c(t,x)u(t,x)
$$
acting on functions given on $\bR^{d+1}$
and elliptic
$$
Mu( x)= a^{ij}( x)u_{x^{i}x^{j}}( x)+
b^{i}( x)u_{x^{i}}( x)+c( x)u( x),
$$
$$
\cM u( x)= \big(a^{ij}(
x)u_{x^{i}}+\hat{b}^{j}(x)
u(x)\big)_{x^{j}}( x)+ b^{i}(
x)u_{x^{i}}( x)+c( x)u( x)
$$
acting on functions given on $\bR^{d}$.
We assume that the coefficients of these operators are measurable
and by magnitude are dominated by a constant $K<\infty$. We also assume
that the matrices $a=(a^{ij})$ are, perhaps, nonsymmetric and satisfy
\begin{equation}
                                                   \label{7.18.1}
a^{ij} \lambda^{i}\lambda^{j}\geq\delta|\lambda|^{2}
\end{equation}
for all $\lambda\in\bR^{d}$ and all possible values of arguments.
Here $\delta>0$ is a fixed constant.

To state our last assumption we set $B_{r}(x)$ to be the open ball 
in $\bR^{d}$ of radius $r$ centered at $x$, $B_{r}=B_{r}(0)$,
$Q_{r}(t,x)=(t,t+r^{2})\times B_{r}(x)$, and $Q_{r}=Q_{r}(0,0)$.
Denote
$$
\osc_{x}(a,Q_{r}(t,x))=
r^{-2}|B_{r}|^{-2}
\int_{t}^{t+r^{2}}\int_{y,z\in B_{r}(x)}|a(s,y)-a(s,z)|\,dydzds,
$$
$$
a^{\#(x)}_{R}=\sup_{(t,x)\in\bR^{d+1}}\sup_{r\leq R}
\osc_{x}(a,Q_{r}(t,x)),\quad a^{\#(x)} =a^{\#(x)}_{\infty}.
$$
This definition is either naturally modified if $a$
is independent of $t$ as in the elliptic operators
or is kept as is. We assume that $a\in VMO_{x}$, that is
$$
\lim_{R\to0}a^{\#(x)}_{R}=0.
$$
For convenience of stating our results we take any 
continuous function $\omega(R)$ on $[0,\infty)$, such that
$\omega(0)=0$ and $a^{\#(x)}_{R}\leq\omega(R)$ for all
$R\in[0,\infty)$. 
Obviously, $a\in VMO_{x}$ if $a$ depends only on $t$. 

By $W^{1}_{p}$ and $W^{2}_{p}$ we denote the usual Sobolev
spaces on $\bR^{d}$.
Also for  $T\in(0,\infty)$
 introduce
$$
\Omega(T)=(0,T)\times\bR^{d}
$$
and, as usual, define $W^{1,2}_{p}(T)$ as the closure of 
the set $C^{1,2}(\Omega(T))$ in the norm
$$
\|u\|_{W^{1,2}_{p}(T)}=\|u\|_{L_{p}(\Omega(T))}
+\|u_{x}\|_{L_{p}(\Omega(T))}+\|u_{xx}\|_{L_{p}(\Omega(T))}
+\|u_{t}\|_{L_{p}(\Omega(T))}.
$$
By $\WO^{1,2}_{p}(T)$ we mean the closure in the same norm
of the subset of $C^{1,2}(\Omega(T))$ consisting of functions
vanishing for $t=T$. Finally,
$$
\cH_{p}(T)=(1-\Delta)^{1/2}W^{1,2}_{p}(T),\quad
\cHO_{p}(T)=(1-\Delta)^{1/2}\WO^{1,2}_{p}(T),
$$
where $\Delta$ is the Laplacian in $x$ variables.
Needless to say   all equations below are
understood in the sense of generalized functions.

Now we fix $T\in(0,\infty)$ and $p\in(1,\infty)$
and state our main results.
\begin{theorem}
                                              \label{theorem 7.20.3}
For any $f\in L_{p}(\Omega(T))$ there exists
a unique $u\in\WO^{1,2}_{p}(T)$ such that
$Lu=f$. Furthermore, there is a constant $N$,
depending only on $d$, $T$, $K$, $\delta$, $p$, and 
the function $\omega$, such that for any $u\in\WO^{1,2}_{p}(T)$
we have
\begin{equation}
                                               \label{7.21.6}
\|u\|_{W^{1,2}_{p}(T)}\leq N\|Lu\|_{L_{p}(\Omega(T))}.
\end{equation}

\end{theorem}
\begin{remark}
                                      \label{remark 7.21.1}
The following is aimed at specialists in 
stochastic processes.
There are solutions of the stochastic differential equation
associated with the operator $L$. We also know that It\^o's
formula is applicable to $u\in W^{1,2}_{p}(T)$ if $p\geq d+1$.
It follows that the solution $u$ is represented as the
expectation of certain integral functional containing $f$.
Such expectations are therefore uniquely defined by
$L$. This leads to weak uniqueness of solutions of
stochastic differential equations with uniformly
nondegenerate bounded diffusion of class $VMO_{x}$ and bounded
measurable drift. More details can be found
in \cite{SV}, where the weak uniqueness is  proved
for equations with uniformly nondegenerate
bounded diffusion,
which is {\em continuous\/} in $x$ uniformly in $t$,
and bounded measurable drift.
\end{remark}
\begin{remark}
                                      \label{remark 7.25.1}
Estimate \eqref{7.21.6} is similar
to interior estimates from \cite{BC}.
However, the space  VMO in \cite{BC} does not include functions
which are independent of $x$ and are measurable in $t$.

\end{remark}
\begin{theorem}
                                      \label{theorem 7.20.4}
Let $f=(f^{1},...,f^{d})$, $g,f^{i}\in L_{p}(\Omega(T))$
for $i=1,...,d$. Then there is a unique
$u\in\cHO_{p}(T)$ such that 
$$
\cL u=\Div f+g.
$$
Furthermore, there is a constant $N$,
depending only on $d$, $T$, $K$, $\delta$, $p$, and 
the function $\omega$, such that 
\begin{equation}
                                            \label{7.25.1}
\|u \|_{L_{p}(\Omega(T))}+
\|u_{x}\|_{L_{p}(\Omega(T))}\leq N
\big(\|f\|_{L_{p}(\Omega(T))}+\|g\|_{L_{p}(\Omega(T))}\big).
\end{equation}
\end{theorem}

\begin{remark}
                                      \label{remark 7.25.2}
Estimate \eqref{7.25.1} is similar
to interior estimates from \cite{By2}.
However, like in \cite{BC}, the space  VMO in \cite{By2}
is defined through approximations by constants and
 does not
include functions which are independent of $x$ and are
measurable in $t$. Also there are no lower order terms in $\cL$
in \cite{By2}.
\end{remark}

\begin{theorem}
                                      \label{theorem 7.20.5}
There exists a constant $\lambda_{0}$, depending only
on  $d$,   $K$, $\delta$, $p$, and 
the function $\omega$, such that, for any $\lambda\geq\lambda_{0}$ 
and  $f\in L_{p}(\bR^{d})$ there exists a unique
$u\in W^{2}_{p}$ satisfying $Mu-\lambda u=f$. 
 
 Furthermore, there is a constant $N$,
depending only on $d$,   $K$, $\delta$, $p$, and 
the function $\omega$, such that for any $u\in W^{2}_{p}$
and $\lambda\geq\lambda_{0}$
we have
\begin{equation}
                                            \label{7.21.5}
\lambda\|u\|_{L_{p}(\bR^{d})}+\|u\|_{W^{ 2}_{p} }\leq
N\|(M-\lambda)u\|_{L_{p}(\bR^{d})}.
\end{equation}
\end{theorem}

\begin{remark}
                                      \label{remark 7.25.3}
Without much stretching the truth one can say that
Theorem \ref{theorem 7.20.5} belongs to
the authors of \cite{CFL1}. The following theorem
is close to some results of \cite{By1}, in which, however,
the lower order terms are not allowed.
\end{remark}

\begin{theorem}
                                        \label{theorem 7.20.6}
There exists a constant $\lambda_{0}$, depending only
on  $d$,   $K$, $\delta$, $p$, and 
the function $\omega$, such that, for any 
$\lambda\geq\lambda_{0}$ 
and  $f=(f^{1},...,f^{d}),g\in L_{p}(\bR^{d})$ there exists a unique
$u\in W^{1}_{p}$ satisfying $\cM u-\lambda u=\Div f+g$. 
 
 Furthermore, there is a constant $N$,
depending only on $\lambda$, $d$,   $K$, $\delta$, $p$, and 
the function $\omega$, such that 
$$
 \|u\|_{W^{ 1}_{p} }\leq
N\big(\|f\|_{L_{p}(\bR^{d})}+\|g\|_{L_{p}(\bR^{d})}\big).
$$
\end{theorem}

\begin{remark}
                                           \label{remark 7.23.1}
As usual in such situations,
from our proofs one can see that instead of the assumption
that
$a\in VMO_{x}$ we are, actually, using that
there   exists $R\in(0,\infty)$ such that
$a^{\#(x)}_{R}\leq\varepsilon$, where $\varepsilon>0$
is a constant depending only on $d,p,\delta,K$.
\end{remark}
\begin{remark}
                                                \label{remark 7.21.2}

Denote  
$$
 u _{Q_{r}(t,x)}=\dashint_{Q_{r}(t,x)} u(s,y)\,dyds,\quad
$$
the average value of a function $u(s,y)$ over $Q_{r}(t,x)$ and
$$
u_{B_{r}(x)}(t)=\dashint_{B_{r}( x)} u(t,y)\,dy 
$$
the average value of a function $u(t,y)$ over $B_{r}(x)$.

Also introduce $\cA$ as the set of $d\times d$ matrix-valued
measurable
functions $\bar{a}=\bar{a}(t)$ depending only on $t$,
satisfying conditions \eqref{7.18.1} and such that
$|\bar{a}^{ij}|\leq K$.
 A standard
fact to recall  is
that for any $\bar{a}\in\cA$
$$
\osc_{x}(a,Q_{r})
\leq2 \dashint_{Q_{r}}|a(s,x)-\bar{a}(s)|\,dxds
$$
and for $\bar{a}(t)= a _{B_{r}}(t)$
$$
\dashint_{Q_{r}}|a(s,x)-\bar{a}(s)|\,dxds\leq\osc_{x}(a,Q_{r}).
$$
This allows one to give obvious equivalent definitions of
$VMO_{x}$.
\end{remark}

\mysection{Auxiliary results}
                                       \label{section 7.22.2}
 
In the  lemmas of this section
\begin{equation}
                                                     \label{7.21.7}
\bar{L}u(t,x)=\bar{a}^{ij}(t)u_{x^{i}x^{j}}(t,x)+u_{t}(t,x),
\end{equation}
where $\bar{a}\in\cA$.

\begin{lemma}
                                                 \label{lemma 7.8.1}
Let $p\in[1,\infty)$,
 $R\in(0,\infty)$,  $u\in
C^{\infty}_{loc}(\bR^{d+1})$,
$$
f=(f^{1},...,f^{d}),\quad  f^{i},g\in L_{p,loc}(\bR^{d+1}),
$$
and  $\bar{L}u=\Div f+g$ in $Q_{R}$. Then
\begin{equation}
                                                     \label{7.8.1}
\int_{Q_{R}}|u(t,x)- u_{Q_{R}}|^{p} \,dxdt
\leq NR^{p} \int_{Q_{R}}(|u_{x }|^{p} +|f|^{p}+R^{p}|g|^{p} )\,dxdt,
\end{equation}
where $N=N(d,K,p)$.
\end{lemma}

Proof. Assume \eqref{7.8.1} is true for $R=1$.
Substitute   $v(t,x)
=u(R^{2}t,Rx)$ in \eqref{7.8.1}  written for $R=1$ and $v$
in place of $u$. Observe that
$$
v_{Q_{1}}=u_{Q_{R}},\quad
\dashint_{Q_{1}}|v(t,x)- v _{Q_{1}}|^{p} \,dxdt
=\dashint_{Q_{R}}|u(t,x)- u _{Q_{R}}|^{p}\,dxdt,
$$
$$
\dashint_{Q_{1}}|v_{x}|^{p}\,dxdt=R^{p}
\dashint_{Q_{R}}|u_{x}|^{p}\,dxdt,
\quad \bar{L}^{R}v(t,x)=R^{2}(\bar{L}u)(R^{2}t,Rx)
$$
$$
=R({\rm div}\,(f(R^{2}t,R\cdot))(x)+R^{2}g(R^{2}t,Rx),
$$
where $\bar{L}^{R}$ is constructed from $\bar{a}(R^{2}t)$.
Then  \eqref{7.8.1} with $R=1$ and $v$ in place of $u$   yields
$$
\dashint_{Q_{R}}|u(t,x)- u_{Q_{R}}|^{p}\,dxdt
\leq NR^{p}\dashint_{Q_{R}}|u_{x}|^{p}\,dxdt
$$
$$
+NR^{p}\dashint_{Q_{R}}|f|^{p}\,dxdt+NR^{2p}
\dashint_{Q_{R}}|g|^{p}\,dxdt.
$$

Hence, we need only prove \eqref{7.8.1} for $R=1$.
In that case take a function $\zeta\in C^{\infty}_{0}(B_{1})$ 
with unit integral. Then by Poincar\'e's inequality,
for any $t\in(0,1)$ and
$$
\bar{u}(t)=\int_{B_{1}}\zeta(y)u(t,y)\,dy.
$$
we have
$$
\int_{B_{1}}|u(t,x)-\bar{u}(t)|^{p}\,dx=
\int_{B_{1}}|\int_{B_{1}}[u(t,x)-u(t,y)]\zeta(y)\,dy|^{p}\,dx
$$
\begin{equation}
                                     \label{7.8.3}
\leq N\int_{B_{1}}\int_{B_{1}}
|u(t,x)-u(t,y)|^{p}\,dxdy
\leq N\int_{B_{1}}|u_{x}(t,x)|^{p}\,dx.
\end{equation}

Observe that for any constant $c$
the left-hand side of
\eqref{7.8.1} is less than a constant times (remember $R=1$)
$$
\int_{Q_{1}}|u(t,x)-c|^{p}\,dxdt\leq2^{p}
\int_{Q_{1}}|u(t,x)-\bar{u}(t)|^{p}\,dxdt+2^{p}
\int_{0}^{1}|\bar{u}(t)-c|^{p}\,dt.
$$
By \eqref{7.8.3} the first term on the right is less
than the right-hand side of \eqref{7.8.1}. To estimate the second term,
take 
$$
c=\int_{0}^{1}\bar{u}(t)\,dt.
$$
Then by Poincar\'e's inequality
$$
\int_{0}^{1}|\bar{u}(t)-c|^{p}\,dt\leq N
\int_{0}^{1}|\int_{B_{1}}\zeta u_{t}\,dx|^{p}\,dt,
$$
where $u_{t}=-(\bar{a}^{ij}u_{x^{i}})_{x^{j}}+{\rm div}\,f+g$.
Integrating by parts with respect to $x$
shows that this term is also less than 
the right-hand side of \eqref{7.8.1}.
The lemma is proved.

\begin{lemma}
                                                 \label{lemma 7.7.2}
There is a constant $N=N(d )$ such that
for any $R>0$ and $u\in C^{\infty}_{loc}(\bR^{d+1})$ we have
\begin{equation}
                                                     \label{7.7.6}
\int_{Q_{R}}|u_{x^{i}}(t,x) 
-(u_{x^{i}})_{Q_{R}}| \,dxdt
\leq NR \int_{Q_{R}}(|u_{xx}| +|u_{t}| )\,dxdt.
\end{equation}
\begin{equation}
                                                     \label{7.7.2}
\int_{Q_{R}}|u(t,x)-u_{Q_{R}}-
x^{i}(u_{x^{i}})_{Q_{R}}| \,dxdt
\leq NR^{2 }\int_{Q_{R}}(|u_{xx}| +|u_{t}| )\,dxdt.
\end{equation}

 \end{lemma}

Proof. To prove \eqref{7.7.6} it suffices to take
$\bar{a}^{ij}=\delta^{ij}$, introduce
$f=\bar{L}u$, note that $\bar{L}(u_{x})=f_{x}$, and apply
 Lemma \ref{lemma 7.8.1} with $p=1$.

To prove \eqref{7.7.2}
set $v(t,x)=u(t,x)- u _{Q_{R}}-
x^{i}(u_{x^{i}})_{Q_{R}}$ and observe that
$$
 v _{Q_{R}}=0,\quad v_{x}=u_{x}-(u_{x})_{Q_{R}}.
$$
Hence for $g:=\bar{L}v$ ($=\bar{L}u$) and $f\equiv0$ by Lemma \ref{lemma 7.8.1}
we find
$$
\int_{Q_{R}}|u(t,x)- u _{Q_{R}}-
x^{i}(u_{x^{i}})_{Q_{R}}| \,dxdt
=\int_{Q_{R}}|v(t,x)-v_{Q_{R}} | \,dxdt
$$
$$
\leq
NR\int_{Q_{R}}(|u_{x}-(u_{x})_{Q_{R}}|+R|u_{t}|+R|u_{xx}|)\,dxdt.
$$
It only remains to use \eqref{7.7.6}. The lemma is proved. 

Define the parabolic boundary of $Q_{r}(t,x)$ by
$$
\partial'Q_{r}(t,x)=\big([t,t+r^{2}]\times\partial B_{r}(x)\big)
\cup\{(t+r^{2},y): y\in  B_{r}(x)\}.
$$
 
\begin{lemma}
                                                 \label{lemma 7.7.3}
Let   $u\in C^{\infty}_{loc}(\bR^{d+1})$
and $\bar{a}$ be infinitely differentiable. Then there
is a unique function $h\in   C^{1,2}(\bar{Q}_{4})$
such that $\bar{L}h=0$ in $Q_{4}$ and   $h=u$
on $\partial'Q_{4}$. Furthermore, $h$ is infinitely differentiable
in $Q_{4}$  and  in $Q_{1}$ we have
\begin{equation}
                                                     \label{7.7.3}
|h_{xx}|+|h_{tx}|+|h_{xxx}| +|h_{txx}| 
\leq N(d,K,\delta) \int_{Q_{4}}(|u_{xx}| +|u_{t}| )\,dxdt.
\end{equation}

\end{lemma}

Proof. The existence,  uniqueness, and the stated properties
 of continuity  $h$ and its derivatives are
classical results.

Therefore, we concentrate on proving \eqref{7.7.3}.
First, notice that by subtracting an appropriate
affine function of $x$ from $u$ and $h$ we reduce the general case to the
one that
\begin{equation}
                                                     \label{7.7.4}
u_{Q_{4}}=(u_{x^{1}})_{Q_{4}}=...=(u_{x^{d}})_{Q_{4}}=0.
\end{equation}
Then, as in the proof of Theorem 8.4.4 of \cite{Kr96} by using Bernstein's
method one proves that for the derivative $D^{\alpha}$ of any order with
respect to
$x$
$$
 \sup_{Q_{1}}|D^{\alpha}h |\leq
N(d,K,\delta,\alpha)\sup_{Q_{2}}|h|.
$$
Since $h_{tx }=-\bar{a}^{ij}h_{ xx^{i}x^{j}}$
and $h_{txx}=-\bar{a}^{ij}h_{xxx^{i}x^{j}}$,
it follows that to prove \eqref{7.7.3} it suffices to prove that

\begin{equation}
                                                     \label{7.7.5}
|h|
\leq N \int_{Q_{4}}(|u_{xx}| +|u_{t}| )\,dxdt
\end{equation}
on $Q_{2}$ under the assumption that \eqref{7.7.4} holds.

Now, take an infinitely differentiable function $\zeta$
on $Q_{4}$ such that it equals 1 near $\partial'Q_{4}$ and zero
inside $Q_{3}$. Without loss of generality assume that
$\bar{a}$ is symmetric.
Then $v=h-\zeta u$ satisfies
$$
\bar{L}v=-\zeta \bar{L}u-2\bar{a}^{ij}\zeta_{x^{i}}u_{x^{j}}-u\bar{L}\zeta
$$
and $v=0$ on $\partial'Q_{4}$. By the maximum principle
$|v|$ is less than the bounded solution $w$ of the Cauchy problem
\begin{equation}
                                                     \label{7.21.1}
\bar{L}w=-(|\zeta \bar{L}u|+2|\bar{a}^{ij}\zeta_{x^{i}}u_{x^{j}}|+
|u\bar{L}\zeta|)I_{Q_{4}}=:f
\end{equation}
in $\{t\leq 16\}$ with zero terminal condition for $t=16$.
This solution is written explicitly as the convolution of 
$f$ and a kernel admitting Gaussian-like
estimates. Since $f$ vanishes inside $Q_{3}$, the convolution
in $Q_{2}$ is estimated by the integral of $f$ over $\bR^{d+1}$.
After that \eqref{7.7.5}
follows from Lemma \ref{lemma 7.7.2}. The lemma is proved.

\begin{lemma}
                                        \label{lemma 7.7.1}
Let $\bar{a}(t)$ be infinitely differentiable.
We assert that
there exists a constant $N=N(d,\delta,K)$ such that
for any $\kappa\geq4$, $r>0$, $u\in
C^{\infty}_{loc}(\bR^{d+1})$, and the solution of $\bar{L}h=0$
in $Q_{\kappa r}$ with boundary condition $h=u$ on
$\partial'Q_{\kappa r}$ we have
\begin{equation}
                                                 \label{7.7.1}
|h_{xx}-( h_{xx})_{Q_{r}}|_{Q_{r}}
\leq N\kappa^{-1} 
(|u_{xx}| +|u_{t}| )_{Q_{\kappa r}} .
\end{equation}

\end{lemma}
 
Proof. Parabolic dilations allow us to only
 concentrate on $r=1$.
The same argument and Lemma \ref{lemma 7.7.3} show that the inequality 
\begin{equation}
                                                     \label{7.9.1}
\kappa|h_{xxx}| +\kappa^{2}|h_{txx}| 
\leq N(d,K,\delta)  (|u_{xx}| +|u_{t}| )
_{Q_{\kappa}} 
\end{equation}
holds in $Q_{\kappa/4}$.
Since $\kappa\geq4$, \eqref{7.9.1} holds in $Q_{1}$ ($=Q_{r}$).
After that it only remains to observe that the 
left-hand side of \eqref{7.7.1} with $r=1$
is less than a constant
times
$$
\sup_{Q_{1}}(|h_{xxx}| +|h_{txx}|)
\leq\kappa^{-1}\sup_{Q_{1}}(\kappa|h_{xxx}| +\kappa^{2}|h_{txx}|)
\leq N\kappa^{-1} (|u_{xx}| +|u_{t}| )_{Q_{\kappa  }}.
$$
The lemma is proved.

\begin{lemma}
                                                 \label{lemma 7.22.1}
Let $q\in(1,\infty)$. Then
there exists a constant
$N=N(q,d,\delta,K)$ such that for any $\kappa\geq4$, $r>0$, $u\in
C^{\infty}_{loc}(\bR^{d+1})$,  we have
\begin{equation}
                                              \label{7.7.01}
|u_{xx}-( u_{xx})_{Q_{r}}| _{Q_{  r}} 
\leq N\kappa^{-1} (|\bar{L}u | +|u_{xx}| )_{Q_{\kappa r}}
+N\kappa^{(d+2)/q}\big(  |\bar{L}u
|^{q} \big)_{Q_{\kappa r}}^{1/q}.
\end{equation}
 \end{lemma}

Proof. We may certainly assume that $\bar{a}$ is 
infinitely differentiable. In that case by Lemma \ref{lemma 7.7.1}
\begin{equation}
                                                     \label{7.22.2}
|h_{xx}-( h_{xx})_{Q_{r}}|_{Q_{r}}
\leq N\kappa^{-1} (|u_{xx}| +|\bar{L}u | )_{Q_{\kappa r}}.
\end{equation}

Furthermore,
$
\bar{L}(u-h)=\bar{L}u
$
in $Q_{\kappa r} $ and $u-h=0$ on $\partial'Q_{\kappa r}$. 
If $\bar{a}$ were constant, then by the standard 
Sobolev space theory we would have
$$
\int_{Q_{\kappa r} }|u_{xx}-h_{xx}|^{q}\,dxdt
\leq N\int_{Q_{\kappa r} }|\bar{L}u|^{q}\,dxdt,
$$
where $N=N(d,\delta,K,q)$.  This estimate is certainly true even
if $\bar{a}$ is not constant. However, we could not find it in
the literature and instead with some reluctance we are going to use Theorem
2.10  of \cite{KK} (see also Remark 2.9 there), which implies that
$$
\int_{Q_{\kappa r} }
(\kappa^{2} r^{2}-|x|^{2})^{q}|u_{xx}-h_{xx}|^{q}\,dxdt
\leq N\int_{Q_{\kappa r} }
(\kappa^{2} r^{2}-|x|^{2})^{q}|\bar{L}u|^{q}\,dxdt,
$$
where $N$ depends only on $\kappa r$, $q$, $d$, $K$,
 and $\delta$.
Observe that $\kappa^{2} r^{2}-|x|^{2}\geq \kappa^{2}
r^{2}/2$ in $Q_{\kappa r/2} $ and $\kappa^{2} r^{2}-|x|^{2}
\leq\kappa^{2} r^{2}$ in $Q_{\kappa r/2} $. It follows that
$$
\int_{Q_{\kappa r/2} } |u_{xx}-h_{xx}|^{q}\,dxdt
\leq N\int_{Q_{\kappa r} }
 |\bar{L}u|^{q}\,dxdt.
$$
Parabolic dilations show that $N$ is independent of $\kappa r$,
and since $\kappa r/2\geq r$ we get
$$
\int_{Q_{ r} }|u_{xx}-h_{xx}|^{q}\,dxdt
\leq N\int_{Q_{\kappa r} }|\bar{L}u|^{q}\,dxdt.
$$
By H\"older's inequality
$$
 |u_{xx}-h_{xx}|_{Q_{r}}
\leq N\kappa^{(d+2)/q}
\big( |\bar{L}u|^{q} \big)_{Q_{\kappa r}}^{1/q} ,
$$
which after being combined with \eqref{7.22.2} shows that there 
is a constant matrix
$\sigma$ ($=(h_{xx})_{Q_{r} }$) such that
$
(|u_{xx}-\sigma|)_{Q_{r}}
$
is less than the right-hand side of \eqref{7.7.01}.
The discussion in the end of Section \ref{section 7.22.1}  shows that
this proves the lemma.

Set
$$
L_{0}u(t,x)=u_{t}(t,x)+a^{ij}(t,x)u_{x^{i}x^{j}}(t,x)
$$
and introduce the (parabolic) maximal and sharp functions of $g$
by
$$
Mg(t,x)=\sup_{r>0}\dashint_{Q_{r}(t,x)}|g(s,y)|\,dyds,
$$
$$
g^{\#}(t,x)=\sup_{r>0}\dashint_{Q_{r}(t,x)}|g(s,y)-
g _{Q_{r}(t,x)}|\,dyds.
$$
 
Here is the main result of this section, in which
all assumptions of Section \ref{section 7.22.1}
are imposed apart from the assumption that $a\in VMO_{x}$.
\begin{theorem}
                                                   \label{theorem 7.9.1}
Let  
$q,\alpha,\beta\in(1,\infty)$,
$\alpha^{-1}+\beta^{-1}=1$, and $R\in(0,\infty)$.
Then 
  there exists a constant $N=N(d,\delta,K,q,\alpha)$ such that
for any $u\in C^{\infty}_{0}(Q_{R})$  we have
\begin{equation}
                                                     \label{7.12.3}
(u_{xx})^{\#}\leq N 
\big[M(|L_{0}u|^{q})\big]^{\mu/q}\big[M |u_{xx}| \big]^{1-\mu} 
+N  \hat{a}^{\mu/(\beta q)}\big[
M(|u_{xx}|^{\alpha q} ) \big]^{1/(\alpha q)},
\end{equation}
on $\bR^{d+1}$,
where 
$ 
\mu=q/(q+d+2),  \hat{a}= a^{\#(x)}_{ R}$. 
\end{theorem}

Proof. 
First, fix  
$\kappa\geq4$,  $r\in(0,\infty)$, and $(t_{0},x_{0})
\in\bR^{d+1}$. Introduce
$$
\bar{a}^{ij}(t)= 
a^{ij}_{B_{\kappa r}(x_{0})}(t) \quad\text{if}\quad 
\kappa r< R,\quad 
\bar{a}^{ij}(t)= 
a^{ij}_{B_{R} }(t) \quad\text{if}\quad 
\kappa r\geq R,
$$
$$
A_{\rho}^{q}= (
|L_{0}u|^{q})_{Q_{\rho}(t_{0},x_{0})},
\quad
B_{\rho}=|u_{xx}|_{Q_{\rho}(t_{0},x_{0})},\quad
 C_{\rho}^{\alpha q}=\big(
 |u_{xx}|^{ \alpha q}\big)_{Q_{\rho}(t_{0},x_{0})},
$$
$$
 A=\sup_{\rho>0}A_{\rho},\quad
B=\sup_{\rho>0}B_{\rho},\quad C=\sup_{\rho>0}C_{\rho}.
$$

By Lemma \ref{lemma 7.22.1} 
$$
|u_{xx}-( u_{xx})_{Q_{r}(t_{0},x_{0})}|_{Q_{r}(t_{0},x_{0})}
 \leq N\kappa^{-1}\ 
(|\bar{L}u |
+|u_{xx}| )_{Q_{\kappa r}(t_{0},x_{0})}
$$
$$
+N\kappa^{(d+2)/q}\big( |\bar{L}u
|^{q} \big)_{Q_{\kappa r}(t_{0},x_{0})}^{1/q}.
$$
By using H\"older's inequality and the fact 
that $\kappa^{-1}\leq1$
we obtain
$$
\kappa^{-1}|\bar{L}u |_{Q_{\kappa r}(t_{0},x_{0})}
\leq\kappa^{(d+2)/q}\big( |\bar{L}u
|^{q} \big)_{Q_{\kappa r}(t_{0},x_{0})}^{1/q},
$$
\begin{equation}
                                              \label{7.22.4}
 |u_{xx}-( u_{xx})_{Q_{r}(t_{0},x_{0})}|_{Q_{r}(t_{0},x_{0})}
\leq N\kappa^{-1}B_{\kappa r} 
+N\kappa^{(d+2)/q}\big( |\bar{L}u
|^{q} \big)_{Q_{\kappa r}(t_{0},x_{0})}^{1/q}.
\end{equation}

Here
$$
\int_{Q_{\kappa r}(t_{0},x_{0})}|\bar{L}u|^{q}\,dxdt
\leq 2^{q}(I+J),
$$
where
$$
I=\int_{Q_{\kappa r}(t_{0},x_{0})}|L_{0}
u|^{q}\,dxdt\leq N(\kappa r)^{d+2}A_{\kappa r}^{q},
$$
$$
J=\int_{Q_{\kappa
r}(t_{0},x_{0})}|(L_{0}-\bar{L})u|^{q}\,dxdt=\int_{Q_{\kappa
r}(t_{0},x_{0})\cap Q_{R}}...
\leq NJ_{1}^{1/\alpha}J_{2}^{1/\beta},
$$
$$
J_{1}=\int_{Q_{\kappa r}(t_{0},x_{0})}
|u_{xx} |^{q\alpha}\,dxdt
\leq N(\kappa r)^{d+2}C_{\kappa r}^{q\alpha},
$$
$$
 J_{2}=\int_{Q_{\kappa r}(t_{0},x_{0})\cap Q_{R}}
|a(t,x)-\bar{a}(t)|^{q\beta}\,dxdt .
$$
If $ \kappa r\geq R$, then we estimate 
$J_{2}$  by the integral over $Q_{R}$, which is less than
$$
  NR^{d+2}\dashint_{  Q_{R}}
|a(t,x)-\bar{a}(t)| \,dxdt 
\leq N (\kappa r) ^{d+2}a^{\#(x)}_{ R}.
$$
In case $  \kappa r<R$  we estimate $J_{2}$ by
$$
 N (\kappa r)  ^{d+2}\dashint_{Q_{ \kappa r}(t_{0},x_{0})}
|a(t,x)-\bar{a}(t)| \,dxdt 
$$
$$
\leq N (\kappa r) ^{d+2}a^{\#(x)}_{\kappa r}
\leq N(\kappa r) ^{d+2}a^{\#(x)}_{ R}.
$$

It follows that
$$
J\leq   N(\kappa r)^{d+2} \hat{a}^{1/\beta}C_{\kappa r}^{q}
$$
and 
$$
\int_{Q_{\kappa r}(t_{0},x_{0})}|\bar{L}u|^{q}\,dxdt
\leq  N(\kappa r)^{d+2}A_{\kappa r}^{q}+
 N(\kappa r)^{d+2} \hat{a}^{1/\beta}C_{\kappa r}^{q}.
$$

Coming back to \eqref{7.22.4} we get
$$
 |u_{xx}-(u_{xx})_{Q_{r}(t_{0},x_{0})}|_{Q_{r}(t_{0},x_{0})}
\leq  N\kappa^{-1}B_{\kappa r}+
N \kappa ^{(d+2)/q}(A_{\kappa r} +
  \hat{a}^{1/(\beta q)}C_{\kappa r})
$$
\begin{equation}
                                        \label{7.31.1}
\leq  N\kappa^{-1}B +
N \kappa ^{(d+2)/q}(A  +
  \hat{a}^{1/(\beta q)}C ).
\end{equation}

So far $\kappa\geq4$ and $r>0$ were fixed.
Now we allow them to vary and observe that \eqref{7.31.1}
is also true for $\kappa\in(0,4)$ since $B $ is
present on the right. After that upon taking supremums  
with respect to $r>0$ and then
minimizing with respect to $\kappa>0$ we
come to
$$
(u_{xx})^{\#}(t_{0},x_{0})\leq N
\big[\hat{a}^{1/(\beta q)}
C +A \big]^{\mu}B^{1-\mu}
$$
$$
\leq N\hat{a}^{\mu/(\beta q)}
C^{\mu}B^{1-\mu }
+NA^{\mu}B^{1-\mu}.
$$
By noting that $B\leq C$ and replacing $B$ with $C$
in the first term on the right we come to what is
  is precisely  
  \eqref{7.12.3} at point $(t_{0},x_{0})$.
The theorem is proved.

Set
\begin{equation}
                                                \label{7.21.2}
L_{p}=L_{p}(\bR^{d+1}).
\end{equation}

\begin{corollary}
                                             \label{corollary 7.15.1}
For any $p\in(1,\infty)$ there exists a constant 
$\varepsilon>0$, depending only on
 $p$, $d$, $K$, and $\delta$, such that   
if $a^{\#(x)}_{R}\leq\varepsilon$ for an $R>0$, then for any
$u\in W^{1,2}_{p}$ we have
\begin{equation}
                                                \label{7.15.1}
\|u_{xx}\|_{L_{p}}\leq N(\|Lu\|_{L_{p}}
+\|u_{x}\|_{L_{p}}+\|u\|_{L_{p}}),
\end{equation}
where
$N=N(R, p,d,K,\delta)$.  
\end{corollary}

Indeed, one the account of the presence of $\|u_{x}\|_{L_{p}}$
and $\|u\|_{L_{p}}$ on the right,
 one may certainly assume that $b\equiv0$ and $c\equiv0$.
The assumption:  $u\in C^{\infty}_{0}(\bR^{d+1})$ also
does not restrict generality.

Next, if $u\in C^{\infty}_{0}(Q_{R})$, then by
 \eqref{7.12.3}, Fefferman-Stein
theorem on sharp functions, and   the Hardy-Littlewood maximal function
theorem
$$
\|u_{xx}\| _{L_{p}}\leq N\|(u_{xx})^{\#}\| _{L_{p}}
\leq N\|[M(|Lu|^{q})]^{1/q}\|_{L_{p}}^{\mu  }
\|u_{xx}\|_{L_{p}}^{ 1-\mu  } 
$$
$$+
N(a^{\#(x)}_{R})^{\mu  /(\beta q)}
\|[M(|u_{xx}|^{\alpha q})]^{1/( \alpha q)}
\| _{L_{p}} 
$$
$$
\leq N\|  Lu  \|_{L_{p}}^{\mu  }
\|u_{xx}\|_{L_{p}}^{ 1-\mu }
+N(a^{\#(x)}_{R})^{\mu  /(\beta q)}\| u_{xx} 
\| _{L_{p}} ,
$$
provided that $p>q\alpha$, that can easily be arranged. 

It follows that
if $a^{\#(x)}_{R}$ is small enough, then 
$$
\|u_{xx}\| _{L_{p}}
\leq N(p,d,K,\delta)\|  Lu  \|_{L_{p}}.
$$  After that 
\eqref{7.15.1} is derived
by a standard procedure using partitions of unity.
We say a little bit more about this procedure in the proof
of Theorem \ref{theorem 7.19.1}.

\mysection{Proof of Theorems \protect\ref{theorem 7.20.3}
and \protect\ref{theorem 7.20.5}}

                                                \label{section 7.22.3}
We suppose that the assumptions of Section \ref{section 7.22.1}
are satisfied
and take a $p\in(1,\infty)$.
We recall notation \eqref{7.21.2} and introduce 
$W^{1,2}_{p}$ as the Sobolev space of functions
$u(t,x)$ on $\bR^{d+1}$ such that $u,u_{x},u_{xx},u_{t}
\in L_{p}$ with natural norm.

\begin{theorem}
                                                   \label{theorem 7.14.1}
There are   constants
$\lambda_{0}$ and
$N$, depending only on $p$, $K$, $\delta$, $d$, and 
$\omega$, such that for any $\lambda\geq\lambda_{0}$
and $u\in W^{1,2}_{p}$ we have  
\begin{equation}
                                                     \label{7.14.3}
\lambda\|u\|_{L_{p}}+\sqrt{\lambda}\|u_{x }\|_{L_{p}}
+\|u_{xx}\|_{L_{p}}+\|u_{t}\|_{L_{p}}\leq N\|(L-\lambda)u\|_{L_{p}}.
\end{equation}
 
Furthermore, for any $\lambda\geq\lambda_{0}$ and
$f\in L_{p}$ there exists a unique 
$u\in W^{1,2}_{p}$ such that $(L-\lambda)u=f$.

\end{theorem}

Proof.   The second assertion is derived from the first one
by the method of continuity. To prove \eqref{7.14.3}
observe that
$$
\|u_{t}\|_{L_{p}}\leq\|Lu \|_{L_{p}}+N\|u_{xx}\|_{L_{p}}
+N\|u_{x}\|_{L_{p}}+\|u\|_{L_{p}},
$$
$$
\|Lu\|_{L_{p}}\leq\|Lu-\lambda u\|_{L_{p}} +\lambda
\| u\|_{L_{p}}.
$$
Therefore,  Corollary
\ref{corollary 7.15.1} shows that  we need only prove that 
for large~$\lambda$
\begin{equation}
                                                     \label{7.15.2}
\lambda\|u\|_{L_{p}}+\sqrt{\lambda}\|u_{x}\|_{L_{p}}  \leq
N\|(L-\lambda)u\|_{L_{p}}.
\end{equation}

We derive \eqref{7.15.2} again from \eqref{7.15.1}  by employing an old
Agmon's idea. Consider the space
$\bR^{d+2}=\{(t,z)=(t,x,y):t,y\in\bR,x\in\bR^{d}\}$ and the function
\begin{equation}
                                                     \label{7.19.3}
\tilde{u}(t,z)=u(t,x)\zeta(y)\cos(\mu y),
\end{equation}
 where
 $\mu=\sqrt{\lambda}$ and $\zeta$ is a $C^{\infty}_{0}(\bR)$-function,
$\zeta\not\equiv0$.
Also introduce the operator
\begin{equation}
                                                     \label{7.19.4}
\tilde{L}u(t,z)=L(t,x)u(t,z)+u_{yy}(t,z).
\end{equation}

Finally, set
$$
\tilde{B}_{r}(z_{0})=\{|z-z_{0}|< r\},\quad
\tilde{Q}_{r}(t_{0},z_{0})=(t_{0},t_{0}+r^{2})\times
\tilde{B}_{r}(z_{0}).
$$

For any $r\in(0,\infty)$, $(t_{0},z_{0})\in\bR^{d+2}$,
and appropriate $\bar{a}(t)$ we have
$$
\int_{ \tilde{Q}_{r}(t_{0},z_{0}) }|a(t,x)-\bar{a}(t)|\,dzdt
$$
$$
\leq\int_{(t_{0}, t_{0}+r^{2})}\int_{
\substack{|x-x_{0}|<r \\
|y-y_{0}|< r}}|a(t,x)-\bar{a}(t)|\,dzdt
$$
\begin{equation}
                                                     \label{7.20.3}
=2r\int_{Q_{r}(t_{0},x_{0}) }|a(t,x)-\bar{a}(t)|\,dxdt
\leq Nr^{d+3}a^{\#(x)}_{R}.
\end{equation}

Since $a\in VMO_{x}$,
it follows that \eqref{7.15.1} holds with $\tilde{u}$, $\tilde{L}$,
and $\bR^{d+2}$ in place of $u$, $L$, and $\bR^{d+1}$, respectively.
Now, observe that
$$
\int_{\bR}|\zeta(y)\sin(\mu y)|^{p}\,dy
$$
is bounded above and away from zero for $\mu\geq1$, so that
$$
\|u_{x }\|^{p}_{L_{p}}
= \big(\int_{\bR}|\zeta(y)\sin(\mu y)|^{p}\,dy\big)^{-1}
\int_{\bR^{d+2}}|u_{x}(t,x)\zeta(y)\sin(\mu y)|^{p}\,dzdt
$$
$$
\leq N \mu^{-p}
\int_{\bR^{d+2}}|u_{x}(t,x)[(\zeta(y)\cos(\mu y))'-
\zeta'(y)\cos(\mu y)]|^{p}\,dzdt
$$
$$
\leq N\mu^{-p}
\int_{\bR^{d+2}}|\tilde{u}_{xy}(t,z)  |^{p}\,dzdt
+N_{1}\mu^{-p}
\int_{\bR^{d+2}}|u_{x}(t,x) 
\zeta'(y)  |^{p}\,dzdt.
$$
The last term can be absorbed by what we started with
if 
$$
N_{1}\mu^{-p}
\int_{\bR }| \zeta'(y)  |^{p}\,dy\leq1/2,
$$
in which case
\begin{equation}
                                                     \label{7.15.4}
\mu\|u_{x }\| _{L_{p} }\leq 
N\|\tilde{u}_{zz}\|_{L{p}(\bR^{d+2})}.
\end{equation}

Similarly,
$$
\|u\|^{p}_{L_{p}}\leq N\mu^{-2p}
\int_{\bR^{d+2}}\big|\tilde{u}_{yy}(t,z)
$$
$$
-u(t,x)[2\mu\zeta'(y)\sin(\mu y)
+\zeta''(y)\cos(\mu y)]\big|^{p}\,dzdt,
$$
$$
\mu^{2} \|u \|_{L_{p} }\leq 
N\|\tilde{u}_{zz}\|_{L{p}(\bR^{d+2})}+N(\mu+1)\|u \|_{L_{p} },
$$
which along with \eqref{7.15.4} yield 
$$
\mu^{2} \|u \| _{L_{p} }+
\mu\|u_{x }\|
_{L_{p} }\leq  N\|\tilde{u}_{zz}\|_{L_{p}(\bR^{d+2})}.
$$

Thus, the left-hand side of \eqref{7.15.2}
is estimated through the left-hand side of
\eqref{7.15.1} written for $\tilde{u}$, $\tilde{L}$,
and $\bR^{d+2}$ in place of $u$, $L$, and $\bR^{d+1}$,
respectively. In turn, the right-hand side of the latter
is easily shown to be less than a constant times
$$
\|\tilde{L}\tilde{u}\|_{L_{p}(\bR^{d+2})}
+\|u_{x}\|_{L_{p}}+(\mu+1)\|u\|_{L_{p}}
$$
$$
\leq N\|Lu-\lambda u\|_{L_{p} }
+\|u_{x}\|_{L_{p}}+(\mu+1)\|u\|_{L_{p}}.
$$
This proves \eqref{7.15.2} and the theorem.

{\bf Proof of Theorem   \ref{theorem 7.20.5}}. As usual, it suffices
to prove the apriori estimate \eqref{7.21.5}. In turn, to do this
it suffices to substitute $v(t,x)=\zeta(t/n)u(x)$, where
$\zeta\in C^{\infty}_{0}(\bR)$, into
\eqref{7.14.3} with $Lv=Mv+v_{t}$, let $n\to\infty$, and observe that
$$
\|v\|^{p}_{L_{p}}=n\alpha\|u\|^{p}_{L_{p}(\bR^{d})},\quad
\|v_{x}\|^{p}_{L_{p}}=n\alpha\|u_{x}\|^{p}_{L_{p}(\bR^{d})},
$$
$$
\|v_{xx}\|^{p}_{L_{p}}=n\alpha\|u_{xx}\|^{p}_{L_{p}(\bR^{d})},\quad
\|(L-\lambda)v\|^{p}_{L_{p}}\leq\|v_{t}\|^{p}_{L_{p}}
$$
$$
+n\alpha\|(M-\lambda)u\|^{p}_{L_{p}(\bR^{d})}
=n^{1-p}\beta\|u\|^{p}_{L_{p}(\bR^{d})}
+n\alpha\|(M-\lambda)u\|^{p}_{L_{p}(\bR^{d})},
$$
where
$$
\alpha=\int_{\bR}|\zeta|^{p}\,dt,\quad
\beta=\int_{\bR}|\zeta'|^{p}\,dt.
$$
The theorem is proved.

{\bf Proof of Theorems  \ref{theorem 7.20.3}}. 
We take $\lambda_{0}$ from Theorem \ref{theorem 7.14.1}.
The method of continuity
and the properties of the heat equation show that
if $g\in L_{p}$ and $g(t,x)=0$ for $t\geq T$, then the solution
$v$ of $(L-\lambda_{0})v=g$ also vanishes for $t\geq T$,
and, therefore, satisfies $(L-\lambda_{0})v=g$ in $\Omega(T)$
with zero condition at $t=T$. We have constructed
a solution from $\WO^{1,2}_{p}(T)$ not of
$Lu=f$ but of $(L-\lambda_{0})v=g$. One gets rid of $\lambda_{0}$
by substitution $u\exp(\lambda_{0}t)=v$.  One   gets
estimate \eqref{7.21.6} from \eqref{7.14.3} by taking
$g=0$ not only for $t\geq T$ but also for $t\leq0$.

It only remains to show uniqueness of solution in
$ \WO^{1,2}_{p}(T)$ of $Lu=0$, which is equivalent
to showing uniqueness for $(L-\lambda_{0})u=0$. We extend
$u$ for $t\geq T$ as zero, obtaining a function
on $(0,\infty)\times\bR^{d}$, which we then extend to negative
$t$ to become an even function of $t$. Call the resulting function
  $\bar{u}$ and denote $f=(L-\lambda_{0})\bar{u}$. Obviously,
  $\bar{u}\in W^{1,2}_{p}$, and since $f=0$ for $t\geq 0$
 the (unique) solution of $ (L-\lambda_{0})v=f$ should 
also vanish for $t\geq0$
as is explained in the beginning of the proof. Hence,
$u=0$ in $\Omega(T)$ and the theorem is proved.

\mysection{Auxiliary results for divergence type equations}
                                             \label{section 7.22.4}

In this section we   discuss some properties
of the operator 
$$
\cL u(t,x)=u_{t}(t,x)+\big(a^{ij}(t,x)u_{x^{i}}(t,x)+
\hat{b}^{j}(t,x)u(t,x)\big)_{x^{j}} 
$$
$$
+b^{i}(t,x)u_{x^{i}}(t,x)+c(t,x)u(t,x). 
$$
All assumptions of Section \ref{section 7.22.1}
are imposed apart from the assumption that $a\in VMO_{x}$.
We   take the operator $\bar{L}$ from \eqref{7.21.7} with
an $\bar{a}\in\cA$.

\begin{lemma}
                                                \label{lemma 7.13.1}
Let $u\in C^{\infty}_{loc}(\bR^{d})$, $f=(f^{1},...
,f^{d})$, $f^{i},g\in L_{2,loc}$ and assume that $\bar{L}u=\Div f+g$
and $\bar{a}$ is infinitely differentiable.
Let $R>0$, $\kappa\geq4$ and let
$h$ be the solution of $\bar{L}h=0$ in $Q_{\kappa R}$ with boundary
condition $h=u$ on $\partial'Q_{\kappa R}$. Then
\begin{equation}
                                                  \label{7.13.1}
\big(|h_{x}-(h_{x})_{Q_{R}}|^{2}\big)_{Q_{R}}
\leq NR^{2}\big(|g|^{2}\big)_{Q_{\kappa R}}
+N\kappa^{-2}\big( |u_{x}|^{2}
+|f|^{2} \big)_{Q_{\kappa R}},
\end{equation}
where $N=N(d,\delta,K)$.
\end{lemma}

Proof. By self-similarity we may assume that $R=1$. Then, by
Lemma \ref{lemma 7.8.1}, applied to $h_{x}$ in place of $u$, 
we see that to prove
\eqref{7.13.1} it suffices to show that
\begin{equation}
                                                  \label{7.13.2}
\big(|h_{xx} |^{2}\big)_{Q_{1}}
\leq N \big(|g|^{2}\big)_{Q_{\kappa  }}
+N\kappa^{-2}\big(|u_{x}|^{2}
+|f|^{2} \big)_{Q_{\kappa  }}.
\end{equation}
The left-hand side of \eqref{7.13.2} will increase if we replace
$Q_{1}$ with $Q_{\kappa/4}$ since $\kappa\geq4$. After that one more
application of parabolic dilations shows that we need only  prove that
\begin{equation}
                                                  \label{7.13.3}
\int_{Q_{1}}|h_{xx} |^{2}\,dxdt\leq
 N \int_{Q_{4}}(|u_{x}|^{2}
+|f|^{2}+|g|^{2})\,dxdt .
\end{equation}
Furthermore, adding a constant  to $u$ results in adding the same
constant to $h$ and does not affect the equation
and \eqref{7.13.3}. Therefore,
we may assume that $u_{Q_{4}}=0$.

If $\zeta\in C^{\infty}_{0}(\bR^{d+1})$ is such that
$\zeta=1$ on $Q_{1}$ and $\zeta=0$ near $\partial'Q_{4}$,
then by observing that
$$
\bar{L}(\zeta h)=2\bar{a}^{ij}\zeta_{x^{i}}h_{x^{j}}+h\bar{L}\zeta
$$
in $\{0<t<16\}$ 
and applying Theorem \ref{theorem 7.20.3}, we see that
the left-hand side of \eqref{7.13.3} is less than  
$$
\int_{0<t<16}|(\zeta h)_{xx}|^{2}\,dxdt\leq N
\int_{0<t<16,\,\zeta\ne0}(|h_{x}|^{2}+
|h|^{2})\,dxdt.
$$
On the account of taking appropriate $\zeta$ we get that
\begin{equation}
                                                  \label{7.18.3}
\int_{Q_{1}}|h_{xx} |^{2}\,dxdt
\leq N\int_{Q_{2}}(|h_{x}|^{2}+
|h|^{2})\,dxdt.
\end{equation}

Then we take a smooth $\eta$ such that $\eta=1$ near $\partial'Q_{4}$
and $\eta=0$ on $Q_{2}$ and observe that  the function
$v= h - u \eta$ vanishes on $\partial'Q_{4}$ and in $Q_{4}$
satisfies
$$
\bar{L}v=-\eta(\Div f+g)- u 
\bar{L}\eta-2\bar{a}^{ij}\eta_{x^{i}}
 u _{x^{j}}.
$$
The usual energy estimate yields
$$
\int_{Q_{4}}(|v_{x}|^{2}+
|v|^{2})\,dxdt\leq N  \int_{Q_{4}}(|u |^{2}+
|f|^{2}+|g|^{2})\,dxdt,
$$
which along with \eqref{7.18.3} lead to
$$
\int_{Q_{1}}|h_{xx} |^{2}\,dxdt
\leq N  \int_{Q_{4}}(|u |^{2}+
|f|^{2}+|g|^{2})\,dxdt.
$$
Finally, Poincar\'e's inequality
(recall that $u_{Q_{4}}=0$) allows us to obtain \eqref{7.13.3}.
The lemma is proved.

\begin{lemma}
                                                 \label{lemma 7.22.2}
Let  
$u\in
C^{\infty}_{loc}(\bR^{d+1})$,$f=(f^{1},...,f^{d})$,
 $f^{i}\in L_{2,loc}$, $\kappa\geq4$, $r>0$.
Assume that $\bar{L}u=\Div f$ in $Q_{\kappa r}$. Then
there exists a constant
$N=N( d,\delta,K)$ such that  
\begin{equation}
                                                 \label{7.22.6}
|u_{x }-( u_{x })_{Q_{r}}|_{Q_{r}}
\leq N\kappa^{-1}
\big(|u_{ x}|^{2} \big)_{Q_{\kappa r}}^{1/2}
+N \kappa^{(d+2)/2}
\big(  |f|^{2}   \big)_{Q_{\kappa r}}^{1/2}.
\end{equation}
 \end{lemma}

Proof. We may assume that $\bar{a}$ is infinitely
differentiable. This follows from the fact that
if $\bar{a}_{n}\in\cA$ are such that $\bar{a}_{n}\to\bar{a}$
(a.e.) as $n\to\infty$ and the operators $\bar{L}_{n}$ are
constructed from $\bar{a}_{n}$, then
 $\bar{L}_{n}u=\Div f_{n}$,
where $f^{i}_{n}=f^{i}+(\bar{a}_{n}^{ij}-\bar{a}^{ij})u_{x^{j}}
\to f^{i}$ in $L_{2}(Q_{\kappa r})$.

Assuming that $\bar{a}$ is infinitely
differentiable we introduce  $h$ as the solution of  
$$
\bar{L} h(t,x):=\bar{a}^{ij}(t)h_{x^{i}x^{j}}(t,x)+h_{t}(t,x)=0
$$
in $Q_{\kappa r} $
with boundary condition $h=u$ on $\partial'Q_{\kappa r} $.
By using Lemma \ref{lemma 7.13.1}
we get (remember $\kappa\geq4$)
$$
\big(|h_{x}-(h_{x})_{Q_{r}}|^{2}\big)_{Q_{r}} 
\leq N\kappa^{-2}\big( |u_{x}|^{2}
+|f|^{2}\big)_{Q_{\kappa r}} 
$$
\begin{equation}
                                                  \label{7.18.6}
 \leq N\kappa^{-2}\big(
|u_{x}|^{2}\big)_{Q_{\kappa r}}
+N\kappa^{d+2}\big(|f|^{2}\big)_{Q_{\kappa r}}.
\end{equation}
 
Furthermore,
 $
\bar{L}(u-h)=\Div f 
 $
in $Q_{\kappa r} $ and $u-h=0$ on $\partial'Q_{\kappa r}$. 
 By the energy
estimate it follows that
$$
\kappa^{-(d+2)}\big(|u_{x }-h_{x }|^{2}\big)_{Q_{ r} }\leq
\big(|u_{x }-h_{x }|^{2}\big)_{Q_{\kappa r} }
\leq 
 N\big(|f|^{2}\big)_{Q_{\kappa r} }.
$$
By combining this with \eqref{7.18.6} and using H\"older's
inequality we see that for the constant vector $\sigma:=(h_{x})_{Q_{r}}$
the expression
$
 |u_{x}-\sigma|_{Q_{r}}
$
is less than the right-hand side of \eqref{7.22.6}. This
proves the lemma.

Now we set
$$
\cL_{0}u(t,x)=u_{t}(t,x)+(a^{ij}(t,x)u_{x^{i}}(t,x))_{x^{j}}
$$
and state the central result of this section.

 \begin{theorem}
                                                   \label{theorem 7.18.1}
Let  
$ \alpha,\beta\in(1,\infty)$,
$\alpha^{-1}+\beta^{-1}=1$, and $R\in(0,\infty)$.
Let  $u\in C^{\infty}_{0}(Q_{R})$, $f=(f^{1},...,f^{d})$,
$f^{i}\in L_{2,loc}$.
Assume that $\cL_{0}u=\Div f$. 
 
Then there exists a constant $N=N(d,\delta,K,\alpha)$ such that
\begin{equation}
                                                     \label{7.18.5}
(u_{x })^{\#} \leq N 
\big[M(|f|^{2})\big]^{\mu/2}\big[M (|u_{x }|^{2} \big]^{(1-\mu)/2} 
+N  \hat{a}^{\mu/(2\beta )}\big[
M(|u_{x }|^{2\alpha}) \big]^{1/(2\alpha  )} 
\end{equation}
on $\bR^{d+1}$, where 
 $
\mu=2/(d+4)$ and $\hat{a}=a^{\#(x)}_{ R}$. 

\end{theorem}

Proof. First, fix  
$\kappa\geq4$, $r\in(0,\infty)$,  and $(t_{0},x_{0})\in\bR^{d+1}$.
 Introduce
$$
A=[M(|f|^{2})(t_{0},x_{0})\big]^{1/2},\quad
B=\big[M (|u_{x}|^{2}) (t_{0},x_{0})\big]^{1/2},
$$
$$ C=\big[
M (|u_{x }|^{2\alpha}) (t_{0},x_{0})\big]^{1/(2\alpha  )}.
$$
Also take $\bar{a}(t) $ as in the proof
of Theorem \ref{theorem 7.9.1} and note that
$$
\bar{L}u=\Div f+\big((\bar{a}^{ij}-a^{ij})u_{x^{i}}\big)_{x^{j}}.
$$

Then by Lemma \ref{lemma 7.22.2}
$$
 |u_{x }-
( u_{x })_{Q_{r}(t_{0},x_{0})}|_{Q_{r}(t_{0},x_{0})}
\leq N\kappa^{-1}B
+N \kappa^{(d+2)/2}(A +D),
$$
where
$$
D^{2}:=\big(|\bar{a}-a|^{2} |u_{x}|^{2}\big
)_{Q_{\kappa r}(t_{0},x_{0})}
\leq
\big(\dashint_{Q_{\kappa r}(t_{0},x_{0})}I_{Q_{R}}|
 \bar{a} -a |^{2\beta }\,dxdt\big)^{1/\beta}
C^{2}.
$$
We estimate the first factor on the right in the same way as in the proof
of Theorem \ref{theorem 7.9.1} and find $D\leq N\hat{a}^{1/(2\beta)}C$,
which leads to
$$
 |u_{x }-
( u_{x })_{Q_{r}(t_{0},x_{0})}| _{Q_{r}(t_{0},x_{0})}
\leq N\kappa^{-1}B
+N\kappa^{(d+2)/2} (A +\hat{a}^{1/(2\beta)}C).
$$
 Having $B$ on the right allows us to assert that this inequality
obtained for $\kappa\geq4$ is actually true for all $\kappa>0$.
Then maximizing with respect to $r>0$ and minimizing
with respect to $\kappa>0$ shows
that
$(u_{x})^{\#}(t_{0},x_{0})$ is less than
$$
N(A^{\mu} +\hat{a}^{\mu/(2\beta)}C^{\mu})B^{1-\mu}.
$$
Observing that $B\leq C$ leads to
\eqref{7.18.5} at $(t_{0},x_{0})$
and proves the theorem.

Similarly to Corollary \ref{corollary 7.15.1} 
we have the following.

\begin{corollary}
                                             \label{corollary 7.20.1}
Let $p\in(2,\infty)$ and $R\in(0,\infty]$.
Then there exist   constants 
$\varepsilon>0$ and $N<\infty$, depending only on
 $p$, $d$, $K$, and $\delta$, such that   
if $a^{\#(x)}_{R}\leq\varepsilon$, then for any
$u\in C^{\infty}_{0}(Q_{R})$ we have
\begin{equation}
                                                \label{7.20.1}
\|u_{x }\|_{L_{p}}\leq N \|f\|_{L_{p}}, 
\end{equation}
 provided that $\cL_{0}u=\Div f$ and $f=(f^{1},...,f^{d})$, $f^{i}\in
L_{p}$.
\end{corollary}

To extend this result to functions not necessarily vanishing
outside $Q_{R}$ we  need to introduce the parameter $\lambda$.

\begin{lemma}
                                             \label{lemma 7.19.1}
Let $p\in(2,\infty)$, $R\in(0,\infty]$,
$f=(f^{1},...,f^{d})$,
$f^{i}, g\in L_{p}(Q_{R})$,
$u\in C^{\infty}_{0}(Q_{R/2}) $, $\lambda\in\bR$, and
$$
\cL u-\lambda u=\Div f+g.
$$
 We assert that
 there exist  constants 
$\varepsilon  \in(0,\infty)$, depending only on
 $p$, $d$, $K$, and $\delta$, and $\lambda_{0},N\in(0,\infty)$,
depending only on the same parameters and $R$, such that   
if $a^{\#(x)}_{R} \leq\varepsilon$, then 
 we have
\begin{equation}
                                                \label{7.19.2}
\sqrt{\lambda}\|u_{x }\|_{L_{p}}+
\lambda \|u \|_{L_{p}}\leq N(\sqrt{\lambda}\|f\|_{L_{p}}
+\|g\|_{L_{p}}),
\end{equation}
provided that   $\lambda\geq\lambda_{0}$.
\end{lemma}

Proof.  First, we observe that
the terms $(\hat{b}^{i}u)_{x^{i}}$ and
$b^{i}u_{x^{i}}$ in $\cL u$ can be included in $\Div f$
and $g$, respectively. Then it is seen that 
without losing generality we may assume that
   $\cL=\cL_{0}$.

In that case  we use the same method as in the proof
of Theorem~\ref{theorem 7.14.1}. 
We take an {\em odd\/} function $\zeta\in C^{\infty}_{0}((-R/2,R/2))$ and
introduce
$\tilde{u}$ and $\tilde{\cL} $ by formulas \eqref{7.19.3} and
\eqref{7.19.4}, of course, taking in the latter $\cL $
in place of $L$. 
Also set $\mu=\sqrt{\lambda}$,
$\tilde{f}^{i}(t,z)=f^{i}(t,x)\zeta(y)\cos(\mu y)$ for
$i=1,...,d$ and
$$
\tilde{f}^{d+1}(t,z)=g(t,x)\zeta_{1}(y)
-2u(t,x)\zeta_{2}(y) +
u(t,x)\zeta_{3}(y),
$$
where
$$
\zeta_{1}(y)=\int_{-\infty}^{y}\zeta(s)\cos(\mu s)\,ds,\quad
\zeta_{3}(y)=\int_{-\infty}^{y}\zeta''(s)\cos(\mu s)\,ds
$$
$$
\zeta_{2}(y)=\mu\int_{-\infty}^{y}\zeta'(s)\sin(\mu s)\,ds
=-\zeta'(y)\cos(\mu y)+\zeta_{3}(y).
$$
Observe that $\zeta_{i}\in C^{\infty}_{0}(\bR)$ since
$\zeta$ is odd and has compact support. Furthermore, as is easy to
check,
$$
\tilde{\cL} \tilde{u}(t,z)
=(\tilde{f}^{1}(t,z))_{x^{1}}+...+
(\tilde{f}^{d}(t,z))_{x^{d}}+(\tilde{f}^{d+1}(t,z))_{y}.
$$

The computations in \eqref{7.20.3} and the fact that
$\tilde{u}$ has support in $(0,R^{2})\times\{|z|<R\}$ convince us that
\eqref{7.20.1} holds for $\tilde{u} $ and $\tilde{f}$  as long as
$a^{\#(x)}$ is small enough.
In other words,
\begin{equation}
                                                \label{7.19.5}
 \|\tilde{u}_{z}\|_{L_{p}(\bR^{d+2})}\leq
N\sum_{i=1}^{d+1}\|\tilde{f}^{i} \|_{ L_{p}(\bR^{d+2})}.
\end{equation}

Since
$$
\kappa_{0}:=\int_{\bR^{d}}|\zeta(y)\sin(\mu y)|^{p}\,dy,\quad
\kappa_{1}:=\int_{\bR^{d}}|\zeta(y)\cos(\mu y)|^{p}\,dy
$$
are bounded  away from zero for   $\mu\geq1$, we get
$$
\|u_{x }\|^{p}_{L_{p} }
=\kappa_{1}^{-1}\int_{\bR^{d+2}}|u_{x}\zeta (y)\cos(\mu y)|^{p}
\,dzdt\leq\kappa_{1}^{-1}\|\tilde{u}_{z}\|^{p}_{L_{p}(\bR^{d+2})},
$$ 
$$
\|u \|^{p}_{L_{p} }=
\kappa_{0}^{-1}\mu^{-1}\int_{\bR^{d+2}}|\tilde{u}_{y}-
u \zeta'(y)\cos(\mu y)|^{p}\,dzdt
$$
$$
\leq N\mu^{-1}(\|\tilde{u}_{z}\|^{p}_{L_{p}(\bR^{d+2})}
+\|u \|^{p}_{L_{p} }).
$$

It follows from here and \eqref{7.19.5} that for $\mu$ large enough
\begin{equation}
                                                \label{7.19.6}
\mu\|u \|^{p}_{L_{p} }+\|u_{x }\|^{p}_{L_{p} }
\leq N\|\tilde{u}_{z}\|^{p}_{L_{p}(\bR^{d+2})}\leq
N\sum_{i=1}^{d+1}\|\tilde{f}^{i} \|_{ L_{p}(\bR^{d+2})}.
\end{equation}

Now we estimate the right-hand side of \eqref{7.19.6}.
Obviously, for
$i=1,...,d$
$$
\|\tilde{f}^{i} \|_{ L_{p}(\bR^{d+2})}\leq N
\|f^{i} \|_{ L_{p} }.
$$

Furthermore,  
$$
\zeta_{1}=\mu^{-1}\big[\zeta (y)\sin(\mu y)-\int_{-\infty}^{y}
\zeta' (s)\sin(\mu s)\,ds\big],
$$
which shows that $\zeta_{1}$ equals $\mu^{-1}$
times a uniformly bounded function with support not wider than that
of $\zeta$. Hence,
$$
\|g\zeta_{1}\|_{ L_{p}(\bR^{d+2})}\leq N\mu^{-1}
\|g \|_{ L_{p} }.
$$
Also $\zeta_{2}$ and $\zeta_{3}$ are  uniformly bounded 
 with support not wider than that
of $\zeta$. Therefore,
$$
\|2u\zeta_{2}+u\zeta_{3}\|_{ L_{p}(\bR^{d+2})}\leq N
\|u\|_{ L_{p} },
$$
$$
\|\tilde{f}^{d+1}\|_{ L_{p}(\bR^{d+2})}\leq N
\mu^{-1}
\|g \|_{ L_{p} }+N\|u\|_{ L_{p} }.
$$
This and \eqref{7.19.6} lead to \eqref{7.19.2}
and the lemma is proved.
\begin{remark}
If $p=2$, then under no restrictions on $a^{\#(x)}$
estimate
  \eqref{7.19.2}  holds for $\lambda$ large, generally,
 or for all
$\lambda>0$ under the additional assumption that $\cL=\cL_{0}$.
This is easily proved by integration by parts.
 \end{remark}

For $n\in\bR$ set
\begin{equation}
                                                \label{7.22.7}
 H^{n}_{p}=(1-\Delta)^{-n/2}L_{p}(\bR^{d}),
\quad
\bH^{n}_{p} = L_{p}(\bR,H^{n}_{p}).
\end{equation} 
 
\begin{theorem}
                                             \label{theorem 7.19.1}
Let $p\in(2,\infty)$,  
$f=(f^{1},...,f^{d})$,
$f^{i},g\in L_{p} $,
$u\in C^{\infty}_{0}(\bR^{d+1}) $, $\lambda\in\bR$,
and
$$
\cL u-\lambda u=\Div f+g.
$$
Take $\varepsilon=\varepsilon(p,d,K,\delta)$
from Lemma \ref{lemma 7.19.1} and assume that
$a^{\#(x)}_{R} \leq\varepsilon$ for an $R\in(0,\infty)$.
Then there exist constants
$\lambda_{0},N\in(0,\infty)$,
depending only on $p,d,K,\delta$, and $R$, such that 
\begin{equation}
                                                \label{7.20.5}
\|u_{t}\|_{\bH^{-1}_{p}}+\sqrt{\lambda}\|u_{x }\|_{L_{p}}+
\lambda \|u \|_{L_{p}}\leq N(\sqrt{\lambda}\|f\|_{L_{p}}
+\|g\|_{L_{p}}) 
\end{equation} 
  whenever    $\lambda\geq\lambda_{0}$.
\end{theorem}

Proof. By the same reasons as before we may assume that $\cL=\cL_{0}$.
Furthermore,
$$
u_{t}=\Div f+g+\lambda u-(a^{ij}u_{x^{i}} )_{x^{i}}
$$
and the operators $(1-\Delta)^{-1/2}$ and
$(1-\Delta)^{-1/2}\partial/(\partial x^{i})$ are
bounded in $L_{p}(\bR^{d})$, which implies that
$$
\|u_{t}\|_{\bH^{-1}_{p}}\leq
N\big(\|f\|_{L_{p}}+\|g\|_{L_{p}}+\lambda\|u\|_{L_{p}}
+\|u_{x}\|_{L_{p}}\big).
$$
It follows that to prove \eqref{7.20.5} it suffices to prove
\begin{equation}
                                                \label{7.20.6}
 \sqrt{\lambda}\|u_{x }\|_{L_{p}}+
\lambda \|u \|_{L_{p}}\leq N(\sqrt{\lambda}\|f\|_{L_{p}}
+\|g\|_{L_{p}}) 
\end{equation} 
and make sure that $\lambda_{0}\geq1$, which is always possible.

Then,
we use partitions of unity. Take a
$\zeta\in C^{\infty}_{0}(Q_{R/2})$ with unit integral,
introduce 
$$
\zeta_{t_{0},x_{0}}(t,x)=\zeta(t-t_{0},x-x_{0}),\quad
u_{t_{0},x_{0}}(t,x)=u(t,x)\zeta_{t_{0},x_{0}}(t,x)
$$
and observe that
$$
(\cL-\lambda) u_{t_{0},x_{0}}=\Div f_{t_{0},x_{0}}
+g_{t_{0},x_{0}},
$$
where
$$
f_{t_{0},x_{0}}^{i}=f^{i}\zeta_{t_{0},x_{0}}
+a^{ji}u(\zeta_{t_{0},x_{0}})_{x^{j}},
$$
$$
 g_{t_{0},x_{0}}=g\zeta_{t_{0},x_{0}}
-f^{i}(\zeta_{t_{0},x_{0}})_{x^{i}}
 +a^{ij}(\zeta_{t_{0},x_{0}})_{x^{j}}u_{x^{i}}
+u(\zeta_{t_{0},x_{0}})_{t}.
$$
It follows from Lemma \ref{lemma 7.19.1} applied to
$Q_{R/2}(t_{0},x_{0})$ in place of $Q_{R/2}$ that
$$
\lambda^{p/2}\|(u_{t_{0},x_{0}})_{x }\|^{p}_{L_{p}}+
\lambda^{p} \|u_{t_{0},x_{0}} \|^{p}_{L_{p}}\leq
N (\lambda^{p/2}\|f_{t_{0},x_{0}}\|^{p}_{L_{p}}
+\|g_{t_{0},x_{0}}\|^{p}_{L_{p}}).
$$
By assuming without losing generality that $\lambda_{0}\geq1$
we estimate the right-hand side by a constant times
$$
 \lambda^{p/2}\|fI_{Q_{R}(t_{0},x_{0})}\|^{p}_{L_{p}}
+\lambda^{p/2}\|uI_{Q_{R}(t_{0},x_{0})}\|^{p}_{L_{p}}
+\|g\zeta_{t_{0},x_{0}}\|^{p}_{L_{p}} 
+ \|u_{x}I_{Q_{R}(t_{0},x_{0})}\|^{p}_{L_{p}} .
$$
On the left
$$
\|u_{x } \zeta_{t_{0},x_{0}}\| _{L_{p}}\leq
\|(u_{t_{0},x_{0}})_{x }\| _{L_{p}}+
\|u(\zeta_{t_{0},x_{0}})_{x}\| _{L_{p}}.
$$
Hence
$$
\lambda^{p/2}\|u_{x } \zeta_{t_{0},x_{0}}\|^{p}_{L_{p}}+
\lambda^{p} \|u\zeta_{t_{0},x_{0}} \|^{p}_{L_{p}}
\leq N\big(
 \lambda^{p/2}\|fI_{Q_{R}(t_{0},x_{0})}\|^{p}_{L_{p}}
$$
$$
+\lambda^{p/2}\|uI_{Q_{R}(t_{0},x_{0})}\|^{p}_{L_{p}}
+\|g\zeta_{t_{0},x_{0}}\|^{p}_{L_{p}} 
+ \|u_{x}I_{Q_{R}(t_{0},x_{0})}\|^{p}_{L_{p}} \big).
$$
After integrating with respect to $(t_{0},x_{0})$ over $\bR^{d+1}$
we conclude
$$
\lambda^{p/2}\|u_{x }  \|^{p}_{L_{p}}+
\lambda^{p} \|u  \|^{p}_{L_{p}}
\leq N\big(
 \lambda^{p/2}\|f \|^{p}_{L_{p}}
$$
$$
+\lambda^{p/2}\|u \|^{p}_{L_{p}}
+\|g \|^{p}_{L_{p}} 
+ \|u_{x} \|^{p}_{L_{p}} \big)
$$
and \eqref{7.19.2} follows. The theorem is proved.

\mysection{Proof of Theorems \protect\ref{theorem 7.20.4}
and \protect\ref{theorem 7.20.6}}
                                             \label{section 7.22.5}

We suppose that the assumptions of Section \ref{section 7.22.1}
are satisfied. 
First, we restate Theorem \ref{theorem 7.19.1} in terms of
appropriate Banach spaces.  We take $\bH^{-1}_{p}$
 from \eqref{7.22.7},
recall that $W^{1,2}_{p}$ is introduced in the beginning
of Section \ref{section 7.22.3} and set
$$
\cH^{1}_{p}=(1-\Delta)^{1/2}W^{1,2}_{p}
$$
with natural norm. It is easy to see that
$C^{\infty}_{0}(\bR^{d+1})$ is dense in $\bH^{n}_{p}$
and $\cH^{1}_{p}$. Furthermore, for 
$u\in C^{\infty}_{0}(\bR^{d+1})$
$$
\|u\|_{\cH^{1}_{p}}=\|(1-\Delta)^{-1/2}u\|_{W^{1,2}_{p}}
\sim\|(1-\Delta)^{-1/2}u\|_{L_{p}}
$$
$$
+\|(1-\Delta)^{-1/2}u_{x}\|_{L_{p}}
+\|(1-\Delta)^{-1/2}u_{t}\|_{L_{p}}
+\|(1-\Delta)^{-1/2}u_{xx}\|_{L_{p}}.
$$
Since the operators $(1-\Delta)^{-1/2}$ and
$(1-\Delta)^{-1/2}(\partial/\partial x^{i})$ are bounded in $L_{p}$,
it follows that
\begin{equation}
                                                       \label{7.22.8}
\|u\|_{\cH^{1}_{p}}\leq N(\|u_{t}\|_{\bH^{-1}_{p}}+
\|u\|_{L_{p}}+\|u_{x}\|_{L_{p}}).
\end{equation}
On the other hand, we know that 
$$
\|u\|_{L_{p}(\bR^{d})}+\|u_{x}\|_{L_{p}(\bR^{d})}
\leq N\|(1-\Delta)^{1/2}u\|_{L_{p}(\bR^{d})}
$$
$$
=N\|(1-\Delta)^{-1/2}(1-\Delta)u\|_{L_{p}(\bR^{d})}
\leq N\|(1-\Delta)^{-1/2}u\|_{L_{p}(\bR^{d})}
$$
$$
+N\|(1-\Delta)^{-1/2}u_{xx}\|_{L_{p}(\bR^{d})},
$$
which shows that the right-hand side of \eqref{7.22.8}
is also dominated by a constant times its left-hand side.
In other words,
$$
\|u\|_{L_{p}}+\|u_{x}\|_{L_{p}}+\|u_{t}\|_{\bH^{-1}_{p}} 
\quad\text{and}\quad\|u\|_{\bH^{1}_{p}}+\|u_{t}\|_{\bH^{-1}_{p}} 
$$
define   equivalent norms in $\cH^{1}_{p}$, which are dominated
by the left-hand side of \eqref{7.20.5} if $\lambda\geq1$.

It turns out that its right-hand side  can be replaced by
a constant times  
\begin{equation}
                                                     \label{8.1.1}
\|\Div f+g\|_{\bH^{-1}_{p}}.
\end{equation}
Indeed, on the one hand, \eqref{8.1.1} is, obviously,
  dominated by the right-hand side of \eqref{7.20.5}.
On the other hand, denote $h=\Div f+g$, $\hat{g}=(1-\Delta)^{-1}h$,
 $\hat{f}^{i}=-\hat{g}_{x^{i}}$. Then
$$
\Div\hat{f}+\hat{g}=-\Delta\hat{g}+\hat{g}=h.
$$
Hence we can replace the right-hand side of \eqref{7.20.5}
with
$$
N(\sqrt{\lambda}\|\hat{f}\|_{L_{p}}+\|\hat{g}\|_{L_{p}}),
$$
which is less than a constant (also depending on $\lambda$)
times $\|h\|_{\bH^{-1}_{p}}$. By the way, the above relations between
$h$, $\hat{g}$, and $\hat{f}$ show that each $h\in\bH^{-1}_{p}$
is written as $\Div\hat{f}+\hat{g}$ 
with $\hat{f}^{i},\hat{g}\in L_{p}$.

It is also worth noting that almost obviously
$\cL$ is a bounded operator
from $\cH^{1}_{p}$ into $\bH^{-1}_{p}$.

This argument allows us to claim that the following result
is a corollary of
Theorem \ref{theorem 7.19.1}.
\begin{theorem}
                                                   \label{theorem 7.22.3}
Let $p\in(2,\infty)$. Then there is a constant
$\lambda_{0}$ depending only on $p$, $d$, $\delta$,
$K$, and $\omega$ such that for any $\lambda\geq\lambda_{0}$
and any $u\in\cH^{1}_{p}$ we have
$$
\|u\|_{\cH^{1}_{p}}\leq N(\lambda,p,d,\delta,K,\omega)
\|(\cL-\lambda)u\|_{\bH^{-1}_{p}}.
$$
\end{theorem}

One derives Theorems  \ref{theorem 7.20.4}
and  \ref{theorem 7.20.6} from the following one
by repeating  the proofs
of   Theorems  \ref{theorem 7.20.3}
and  \ref{theorem 7.20.5} almost word for word.

\begin{theorem}
                                                  \label{theorem 7.18.2}
Let $p\in(1,\infty)$. Then there exists
  $\lambda_{0}=\lambda_{0}(p,d,K,\delta,\omega)$
such that,    for any
$u\in\cH^{1}_{p} $ and $\lambda \geq\lambda_{0}$,   we have
\begin{equation}
                                                       \label{7.18.8}
\|u_{t}\|_{\bH^{-1}_{p} }+\|u\|_{\bH^{1}_{p} }\leq 
N(p,d,K,\delta,R,\lambda)\|(\cL -\lambda)
u\|_{\bH^{-1}_{p} }.
\end{equation}
 
Furthermore, for each  $\lambda\geq\lambda_{0}$ and
  $f\in \bH^{-1}_{p} $
there is a unique
$u\in\cH^{1}_{p} $ such that   $(\cL-\lambda) u=f$.

\end{theorem}

Proof. The second assertion is a standard consequence of 
the first one and
the method of continuity. 

 Estimate \eqref{7.18.8} is stated in Theorem \ref{theorem 7.22.3}
for $p>2$. 

We consider $p\in(1,2)$ by using duality. Set $q=p/(p-1)$,
take $\lambda_{0}$ corresponding to $q$ and the operator
$\cL^*$ formally adjoint to $\cL$ and take 
$\lambda\geq\lambda_{0}$. The reader should not  
be uncomfortable with the fact that the derivative in time
enters $\cL^*$ with a negative sign unlike $\cL$. Our results are
applicable to
such operators as well, which is seen after changing variables
$t\to-t$.

By the above
for any $h\in\bH^{-1}_{q}$ we can find $v\in\cH^{1}_{q}$
such that
$$
(\cL^{*}-\lambda)v=h,\quad\|v\|_{\bH^{1}_{q}}\leq N\|h\|_{\bH^{-1}_{q}}.
$$
For $u\in\cH^{1}_{p}$ write
$$
|(u,h)|=|((\cL-\lambda)u,v)|\leq\|(\cL-\lambda)u\|_{\bH^{-1}_{p}}
\|v\|_{\bH^{1}_{q}}\leq
N\|(\cL-\lambda)u\|_{\bH^{-1}_{p}}\|h\|_{\bH^{-1}_{q}}.
$$
Since $h$ was arbitrary, it follows that
$$
\|u\|_{\bH^{1}_{p} }\leq 
N \|(\cL -\lambda)
u\|_{\bH^{-1}_{p} }.
$$
This estimate and the formula
$$
u_{t}=(\cL-\lambda)u+\lambda u-(a^{ij}u_{x^{i}}+\hat{b}^{i}u)_{x^{j}}
-b^{i}u_{x^{i}}
$$
allow us to get the remaining part of \eqref{7.18.8},
which is thus proved for $p\in(1,\infty)$, $p\ne2$. 
Once the resolvent operator is constructed for $p\in(1,\infty)$,
$p\ne2$,
the case
$p=2$ is covered by interpolation. The theorem is proved.

\end{document}